\documentclass[12pt]{article}

\usepackage{amsfonts}

\newtheorem{thm}{Theorem}[section]

\newtheorem{cor}[thm]{Corollary}

\newtheorem{prop}[thm]{Proposition}
\newtheorem{conjecture}{Conjecture}

\newcommand{\be}{\begin{equation}}
\newcommand{\ee}{\end{equation}}
\newcommand{\openbox}{\leavevmode
  \hbox to8pt{\hfil\vrule\vbox to6pt{\hrule width6pt\vfil\hrule}\vrule}}

\newcommand{\ve}[1]{\mathbf{#1}}
\newcommand{\qed}{\hbox to5pt{ } \hfill \openbox\bigskip\medskip}

\newcommand{\Zp}{\mathbb Z _p}

\newcommand{\cF}{\mbox{$\cal F$}}
\newcommand{\cG}{\mbox{$\cal G$}}

\newcommand{\cB}{\mbox{$\cal B$}}

\newcommand{\Sf}{\mathbb S}

\newcommand{\Z}{\mathbb Z}

\newcommand{\R}{\mathbb R}

\newcommand{\F}{\mathbb F}

\title{A new upper bound for the size of $s$-distance sets in boxes}
\author{G\'abor Heged\"{u}s
\\{\normalsize  \'Obuda University}
\\{\normalsize B\'ecsi \'ut 96, Budapest, Hungary, H-1037}
\\{\normalsize hegedus.gabor@nik.uni-obuda.hu}
}

\begin{document}

\maketitle
\begin{abstract}
Let $q,d\geq 2$ be integers. Define
$$
J(q,d):=\frac 1q \Big( \min_{0<x<1} \frac{1-x^q}{1-x} x^{-\frac{q-1}{d}}\Big).
$$
Let  $\mbox{$\cal G$}\subseteq {\mathbb R}^n$ be an arbitrary subset. We denote by $d(\mbox{$\cal G$})$ the set of (non-zero) distances among points of $\mbox{$\cal G$}$:
$$
d(\mbox{$\cal G$}):=\{d( p_1, p_2):~  p_1, p_2\in \mbox{$\cal G$}, p_1\ne  p_2\}.
$$

Our main result is a new upper bound for the size of $s$-distance sets in boxes. More concretely, 
let  $A_i\subseteq \mathbb R$,  $|A_i|=q\geq 2$ be subsets for each $1\leq i\leq n$. Consider the box $\mbox{$\cal B$}:=\prod_{i=1}^n A_i\subseteq  {\mathbb R}^n$. Suppose that $\mbox{$\cal G$}\subseteq \mbox{$\cal B$}$ is a set such that $|d(\mbox{$\cal G$})|\leq s$. Let $d:=\frac{n(q-1)}{s}$. Then
$$
|\mbox{$\cal G$}|\leq 2(qJ(q,d))^n.
$$

We use  Tao's slice rank bounding method in our proof.
\end{abstract}

\medskip

\section{Introduction}

In this manuscript we give upper bounds for $s$-distance sets in the direct product of finite set of points in the Euclidean space.   Two independent research directions  motivated our results. 

Our first motivation comes from the following  question of Erd\H{o}s:  Given $n$ points in the plane, what is the smallest number of distinct distances they can determine?

Erd\H{o}s proved in  \cite{E} the following result.

\begin{thm} \label{Erdos}
The minimum number of $f(n)$ of distances determined  by $n$ points of the plane satisfies the inequalities 
$$
(n-\frac 34)^{1/2}- \frac 12\leq f(n)\leq \frac{cn}{(\log n)^{1/2}}.
$$
\end{thm}
Erd\H{o}s conjectured that the square grid is essentially the extremal example, consequently the upper bound for the function $f(n)$ is sharp.   

Our second motivation comes from the investigation of different types of  $s$-distance sets  in the $d$-dimensional Euclidean space.

First we introduce some notation. Let ${\R}^d$ denote the $d$-dimensional Euclidean space. Let  $\cG\subseteq {\R}^n$ be an arbitrary set. We denote by $d(\cG)$ the set of (non-zero) distances among points of $\cG$:
$$
d(\cG):=\{d(\ve p_1,\ve p_2):~ \ve p_1,\ve p_2\in \cG,\ve p_1\ne \ve p_2\}.
$$

Let $(\ve x,\ve y)$ stand for  the standard scalar product. Let  $\cG\subseteq {\R}^n$ be an arbitrary set. We denote by $s(\cG)$ the set of scalar products between the dintinct points of $\cG$. 
$$
s(\cG):=\{(\ve p_1,\ve p_2):~ \ve p_1,\ve p_2\in \cG,\ve p_1\ne \ve p_2\}.
$$

Bannai, Bannai and  Stanton proved the following result in \cite{BBS} Theorem 1.

\begin{thm} \label{BBSupper}
Suppose that $\cF\subseteq {\R}^n$ is a set and that $|d(\cF)|\leq s$. Then
$$
|\cF|\leq {n+s\choose s}.
$$
\end{thm}

Delsarte, Goethals and  Seidel investigated the spherical $s$-distance sets. They proved  the following Theorem in \cite{DGS}.
\begin{thm} \label{DGGupper}
Suppose that $\cF\subseteq {\Sf}^{n-1}$ is a set and that $|d(\cF)|\leq s$.
Then
$$
|\cF|\leq {n+s-1 \choose s}+ {n+s-2 \choose s-1}.
$$
\end{thm}

We state here our main results.
\begin{thm} \label{main}
Let  $A_i\subseteq \R$,  $|A_i|=q\geq 2$ for each $1\leq i\leq n$. Consider the box $\cB:=\prod_{i=1}^n A_i\subseteq  {\R}^n$. Suppose that $\cF\subseteq \cB$ is a set such that $|d(\cF)|\leq s$.
Then
$$
|P|\leq 2\cdot |\{x_1^{\alpha_1}\cdot\ldots \cdot x_n^{\alpha_n}:~ 0\leq \alpha_i \leq q-1 \mbox{ for each } i,\ \sum_i \alpha_i \leq s  \}|.
$$
\end{thm}

We use  Tao's slice rank bounding method in our proof (see  the blog \cite{T}). Tao developed this method as a proof technique to prove Ellenberg and Gijswijt's breakthrough about the upper bounds for the size of subsets $A$ in $({\Zp})^n$ without three-term arithmetic progressions (see \cite{EG}). 

Let $t,d\geq 2$ be integers. Define
$$
J(t,d):=\frac 1t \Big( \min_{0<x<1} \frac{1-x^t}{1-x} x^{-\frac{t-1}{d}}\Big).
$$
Define $J(q):=J(q,3)$ for each $q>1$.

This $J(q)$ constant appeared in the proof of the Ellenberg and Gijswijt's bound for the size of three-term progression-free sets (see \cite{EG}). It was proved in \cite{BCCGU} Proposition 4.12 that $J(q)$ is a decreasing function of $q$ and
$$
\lim_{q\to \infty} J(q)=\inf_{z>3} \frac{z-z^{-2}}{3\log(z)}=0.8414\ldots .
$$ 
We can verify easily that $J(3)=0.9184$, consequently $J(q)$ lies in the range 
$$
0.8414\leq J(q)\leq 0.9184
$$
for each $q\geq 3$.
The following Corollary gives a more concrete upper bound for the size of  $s$-distance sets in boxes. 
\begin{cor} \label{cor}
Let  $A_i\subseteq \R$,  $|A_i|=q\geq 2$ for each $1\leq i\leq n$. Consider the box $\cB:=\prod_{i=1}^n A_i\subseteq  {\R}^n$. Suppose that $\cF\subseteq \cB$ is a set such that $|d(\cF)|\leq s$.
Let $d:=\frac{n(q-1)}{s}$. Then
$$
|\cF|\leq 2(qJ(q,d))^n.
$$
\end{cor}

Deza and Frankl proved the following statement in \cite{DF} Theorem 4. 
\begin{thm} \label{DFrank}
Suppose that $\cF\subseteq {\R}^n$ is an arbitrary  set such that $|s(\cF)|\leq s$. Then
$$
|\cF|\leq {n+s\choose s}.
$$
\end{thm}

It is easy to prove the following result using a slight modification of Tao's slice rank bounding method.
\begin{thm} \label{DFrankbox}
Let  $A_i\subseteq \R$,  $|A_i|=q\geq 2$ for each $1\leq i\leq n$. Consider the box $\cB:=\prod_{i=1}^n A_i\subseteq  {\R}^n$. Suppose that $\cF\subseteq \cB$ is an arbitrary set which satisfies the  following properties:
\begin{itemize}
\item[(i)] $(\ve f,\ve f)\notin s(\cF)$ for each $\ve f\in \cF$;
\item[(ii)] $|s(\cF)|\leq s$.
\end{itemize}
Then 
$$
|\cF|\leq |\{x_1^{\alpha_1}\cdot\ldots \cdot x_n^{\alpha_n}:~ 0\leq \alpha_i \leq q-1 \mbox{ for each } i,\ \sum_i \alpha_i \leq s  \}|.
$$
\end{thm}
In Section 2 we present our proofs.

\section{Proofs}

The following simple statement was proved in \cite{CLP}  Lemma 1.
\begin{prop} \label{CLP}
Suppose that $n \geq 1$ and $d \geq 0$ are integers, $P$ is a multilinear polynomial in $n$ variables of total degree at most $d$ over a field $\F$, and $\cF \subseteq {\F}^n$ is a subset with 
$$
|\cF|>2\sum_{i=0}^{d/2} {n\choose i}.
$$
If $P(\ve a- \ve b)=0$ for all $\ve a, \ve b\in \cF$, $\ve a\neq \ve b$, then $P(\ve 0)=0$.
\end{prop}

Tao's slice rank bounding method (see the blog \cite{T} and an other proof in \cite{BCCGU} Section 4) gives easily the following generalization of Proposition \ref{CLP}, which is a special case of \cite{H} Corollary 1.3.

\begin{thm} \label{maincor}
Let $\F$ be an arbitrary field. Let $A_i\subseteq \F$ be fixed subsets   such that $|A_1|= \ldots =|A_n|=t>0$. Let $\cF\subseteq \prod_{i=1}^n A_i$ be a finite subset. 
Suppose that there exists a polynomial 
$$
P(x_{1}, \ldots ,x_{n},y_{1}, \ldots, y_{n}, )\in \F[x_{1}, \ldots x_{n},y_{1}, \ldots y_{n} ]
$$
satisfying the following conditions:
\begin{itemize}
\item[(i)] $P(\ve a,\ve a)\ne 0$ for each  $\ve a\in \cF$;
\item[(ii)] if $\ve a,  \ve b\in \cF$, $\ve a\ne \ve b$ are arbitrary vectors, then $P(\ve a,\ve b)=0$.
\end{itemize}
Then 
$$           
|\cF| \leq 2|\{x_1^{\alpha_1}\ldots x_n^{\alpha_n}:~ 0\leq {\alpha_i\leq t-1} \mbox{ for each } 1\leq i\leq n  \mbox{ and}
$$
\be \label{upper}
\sum_{j=1}^n \alpha_j\leq \frac{\mbox{deg}(P)}{2} \} |.
\ee
\end{thm}

As an easy consequence, we proved the following result in \cite{H} Corollary 1.5.

\begin{thm} \label{maincor2}
Let $\F$ be an arbitrary field. Let $A_i\subseteq \F$ be fixed subsets  such that $|A_1|= \ldots =|A_n|=t>0$. Let $\cF\subseteq \prod_{i=1}^n A_i$ be a finite subset. 
Suppose that there exists a polynomial 
$$
P(x_{1}, \ldots ,x_{n},y_{1}, \ldots, y_{n}, )\in \F[x_{1}, \ldots x_{n},y_{1}, \ldots y_{n} ]
$$
satisfying the following conditions:
\begin{itemize}
\item[(i)] $P(\ve a,\ve a)\ne 0$ for each  $\ve a\in \cF$;
\item[(ii)] if $\ve a,  \ve b\in \cF$, $\ve a\ne \ve b$ are arbitrary vectors, then $P(\ve a,\ve b)=0$.
\end{itemize}
Let $d:=\frac{2n(t-1)}{deg(P)}$.  Then 
$$   
|\cF| \leq 2(tJ(t,d))^n.
$$
\end{thm}

{\bf Proof of Theorem \ref{main}:}

Consider the set $d(\cF)=\{d_1, \ldots , d_s\}$. Clearly $d_i\ne 0$ for each $i$.

Define the polynomial
$$
P(x_1, \ldots ,x_n, y_1, \ldots ,y_n):= \prod_{i=1}^s \Big(\sum_{j=1}^n (x_j-y_j)^2-d_i^2 \Big)\in \R[x_1, \ldots ,x_n, y_1, \ldots ,y_n].
$$
Clearly $deg(P)=2s$.
Then
$$
P(\ve a,\ve a)=\prod_{i=1}^s (-d_i)^2\ne 0
$$
for each $\ve a\in \cF$
On the other hand  if $\ve a, \ve b\in \cF$, $\ve a\neq \ve b$, then 
$$
P(\ve a,\ve b)=\prod_{i=1}^s \Big(\sum_{j=1}^n (a_j-b_j)^2-d_i^2\Big)=\prod_{i=1}^s \Big(d(\ve a,\ve b)^2-d_i^2\Big).
$$

But there exists an $i$ such that $d(\ve a,\ve b)=d_i$, because $d(\cF)=\{d_1, \ldots , d_s\}$. Hence $d(\ve a,\ve b)^2=d_i^2$, so $P(\ve a,\ve b)=0$.

Finally we can apply Theorem \ref{maincor} with the choices $\F=\R$ and  $t=q$. \qed

\section{Concluding remarks}

We conjecture that the following upper bound is sharp. 

\begin{conjecture} \label{upperb}
Let  $A_i\subseteq \R$,  $|A_i|=q\geq 2$ for each $1\leq i\leq n$. Consider the box $\cB:=\prod_{i=1}^n A_i\subseteq  {\R}^n$. Suppose that $\cF\subseteq \cB$ is a set with $|d(\cF)|\leq s$.
Then
$$
|\cF|\leq |\{x_1^{\alpha_1}\cdot\ldots \cdot x_n^{\alpha_n}:~ 0\leq \alpha_i \leq q-1 \mbox{ for each } i,\ \sum_i \alpha_i \leq s  \}|.
$$
\end{conjecture}
As a lower bound, we give the following construction in the case $q=2$. If $A\subseteq [n]$, then denote by $\ve v_A$ the characteristic vector of $A$. Consider the set
$$
\cF:=\{\ve v_A:~ A\subseteq [n], |A|= s\}\subseteq \{0,1\}^n\subseteq {\R}^{n}.
$$ 
Then $\cF$ is an $s$-distance set with $|\cF|={n \choose s}$.


\end{document}